\documentclass[11pt]{amsart}

\usepackage{amssymb} 
\usepackage[mathscr]{eucal} 
\usepackage{xypic}   
\usepackage{amsmath} 
\usepackage{epsfig}
\usepackage{amscd}

\newtheorem{TEO}{Theorem}[section]
\newtheorem{PROP}[TEO]{Proposition}
\newtheorem{LEM}[TEO]{Lemma}

\newtheorem{REM}[TEO]{Remark}

\newcommand\dual{\mathrel{\raise3pt\hbox{$\underline{\mathrm{\thinspace d
\thinspace}}$}}}

\newcommand\proj{\mathbb P}
\newcommand\Z{\mathbb Z}

\def\Z{{\mathbb Z}}

\def\Z{{\mathbb Z}}

\begin{document}

\title{Some results on the second Gaussian map for curves}
\author{Elisabetta Colombo\\ Paola Frediani}
\date{}
%
%


\footnote{Elisabetta Colombo, Dipartimento di Matematica,
Universit\`a di Milano, via Saldini 50,
     I-20133, Milano, Italy. e-mail: elisabetta.colombo@mat.unimi.it

Paola Frediani, Dipartimento di Matematica, Universit\`a di Pavia,
via Ferrata 1, I-27100 Pavia, Italy. e-mail:
paola.frediani@unipv.it

The present research took place in the framework of the
PRIN 2005 of MIUR: "Spazi dei moduli e teoria di Lie" and PRIN 2006 of MIUR "Geometry on algebraic varieties". 

 }
\maketitle

\begin{abstract}
We study  the second Gaussian map  for a curve of genus $g$, in
relation with the second fundamental form of the period map
$j:M_g\rightarrow A_g$. We exhibit a class of infinitely many curves with surjective second Gaussian map.
We compute its rank on the hyperelliptic
and trigonal locus and show global generation of its image in
$H^0(X,4K_X)$ for $X$ not hyperelliptic nor trigonal.
\end{abstract}

\section{Introduction}

The first Gaussian map for the canonical series has been
intensively studied, and it has been shown that
 for a general curve of genus different from 9 and $\leq 10$  it is injective, while for genus $\geq 10$,
 different from 11, it is surjective (\cite{chm}, \cite{voi}, \cite{cm}). In \cite{wahl3} it is proven that
 if a curve lies on a K3 surface, the first Gaussian map can't be surjective and it is known (\cite{mm}) that
 the general curve of genus 11 lies on a K3 surface.

 In this paper we study some properties of the second Gaussian map
$\mu_2:I_2(K_X)\rightarrow H^0(X,4K_X)$. Our geometrical
motivation comes from its relation with the curvature of the
moduli space  ${M_g}$ of curves of genus $g$ with the Siegel
metric induced by the period map $j:{ M_g}\rightarrow {A_g}$, that
we started to analyze in \cite{cfcl}.

There the curvature is computed using the formula for the
associated second fundamental form given in \cite{cpt}. In
particular in \cite{cpt} it is proven that the second fundamental
form lifts the second Gaussian map $\mu_2$ as stated in an
unpublished paper of Green-Griffiths (cf. \cite{green}).

In \cite{cfcl}, (3.8) we give a formula for the holomorphic
sectional curvature of ${M_g}$ along the
 a Schiffer variation $\xi_P$, for $P$ a point  on the curve
 $X$, in terms of the holomorphic sectional curvature of ${A_g}$ and the second Gaussian map.

The relation of the second Gaussian map with curvature properties
of $M_g$ in $A_g$ suggests that its rank could give information on
the geometry of $M_g$. Note that surjectivity can be expected for
a general curve of genus at least 18. Recall that $M_g$ is
uniruled for $g \leq 15$, it has Kodaira dimension at least 2 for
$g=23$, and it is of general type for all other values of  $g \geq
22$.

In this direction in this paper we exhibit infinitely many
examples of curves lying on the product of two curves with
surjective second Gaussian map. Other examples of curves whose
second Gaussian map is surjective were given by in \cite{bafo} for
complete intersections. Both classes of examples generalize
constructions given by Wahl for the first Gaussian map in
\cite{wahl1}, \cite{wahl3}.

We are also able to determine the rank of $\mu_2$ on the
hyperelliptic and trigonal loci.
More precisely for any hyperelliptic curve of genus $g\geq 3$ we
show that
 $rank(\mu_2)=2g-5$ and its image has the Weierstrass points as base points.
 For any trigonal (non hyperelliptic) curve of genus $g\geq 8$ we show that
 $rank(\mu_2)=4g-18$ and its image has the ramification points of the $g^1_3$ as base points.

Finally we prove that for any non-hyperelliptic, non-trigonal
curve of genus $g\geq
 5$ the image of $\mu_2$ has no base points.

In \cite{cfcl} we apply these results to the holomorphic sectional
curvature of ${M_g}$. In particular along a Schiffer variation
$\xi_P$ the holomorphic sectional curvature $H(\xi_P)$ of ${M_g}$
is strictly smaller than the holomorphic sectional curvature of ${
A_g}$ for a non-trigonal, nor hyperelliptic curve
(\cite{cfcl},(4.5)). Instead, if $P$ is either a Weierstrass point
of a hyperelliptic curve or a ramification point of the $g^1_3$ on
a trigonal curve, the holomorphic sectional curvature $H(\xi_P)$
is equal to the holomorphic sectional curvature of ${ A_g}$, which
equals $-1$ (\cite{cfcl}, (4.5), (5.3)).

The computations are based on the observation that,
 for a quadric $Q$ of rank at most 4,
 $\mu_2(Q)$ can be written as the product of the first Gaussian
 maps associated to sections of the two adjoint line bundles $L$ and
 $K\otimes L^{-1}$ which define the quadric.
 As a first straightforward consequence, we show that for any curve the rank of $\mu_2$
 is greater or equal to $g-3$.

To study the trigonal case we use the related results on the first Gaussian map for trigonal curves of \cite{md} and \cite{br}.
A crucial step to determine the rank of the first and hence
the second Gaussian map on the trigonal locus is the observation
that a trigonal curve lies on a rational normal scroll.
A natural question is to understand whether
restrictions on the rank of $\mu_2$ can be obtained if a curve
lies on a special surface, as it happens for the first Gaussian
map and $K3$ surfaces, or if it occurs in a non trivial linear
system of a surface (cf. \cite{ser} for related results). We
intend to continue our investigations on the rank properties of
$\mu_2$ for general curves and for curves on surfaces in the near
future.

 The paper is organized as follows:

 In Section 2 we study the second Gaussian map on quadrics of rank less or equal to 4, in terms of the first Gaussian maps associated to sections of the two adjoint line bundles $L$ and
 $K\otimes L^{-1}$ which define the quadric.

In Section 3 we show a class of infinitely many examples of curves with surjective second Gaussian map.

 In Section 4 we determine the rank of $\mu_2$ for hyperelliptic and trigonal curves.

In Section 5 we prove injectivity of $\mu_2$ for general curve of genus at most six by specialization on trigonal curves
and on the smooth plane quintic.

 In Section 6 we study the global generation of the image of $\mu_2$.

{\bf Acknowledgments.} The authors thank Gilberto Bini and Pietro
Pirola for several fruitful suggestions and discussions on the
subject.

\section{The Second Gaussian map}

We first recall the definition of the Gaussian maps
(cf. \cite{wahl2}). Let $X$ be a smooth projective curve, $S := X
\times X$, $\Delta \subset S$ be the diagonal. Let $L$ be a line
bundle on $X$ and $L_S := p_1^*(L) \otimes p_2^*(L)$, where $p_i :
S \rightarrow X$ are the natural projections. Consider the
restriction map
$$\tilde{\mu}_{n,L}: H^0(S, L_S(-n\Delta)) \rightarrow H^0(\Delta ,
L_S(-n\Delta)_{| \Delta}).$$
Notice that since ${\mathcal O}(\Delta)_{|\Delta} \cong T_X$, we have
$$H^0(\Delta ,L_S(-n\Delta)_{| \Delta}) \cong H^0(X, 2L \otimes nK_X).$$
In the case $L = K_X$, $I_2(K_X) \subset H^0(S, K_S(-2\Delta))$,
so we can define the second Gaussian map
$$\mu_2:  I_2(K_X) \rightarrow H^0(X, 4K_X),$$
as the restriction $\tilde{\mu}_{2,K|I_2(K_X)}$.

As above we fix a basis $\{\omega_i\}$ of $H^0(K_X)$. In local coordinates
$\omega_i = f_i(z) dz$. Let $Q \in I_2(K_X)$, $Q = \sum_{i,j} a_{ij} \omega_i \otimes
\omega_j$, recall that $\sum_{i,j} a_{ij} f_if_j \equiv 0$, and since
$a_{i,j}$ are symmetric, we also have $\sum_{i,j} a_{ij} f'_if_j \equiv 0$.
The local expression of $\mu_2(Q)$ is
\begin{equation}
\mu_2(Q) = \sum_{i,j} a_{ij} f''_if_j (dz)^4= - \sum_{i,j}a_{ij} f'_i f'_j(dz)^4.
\end{equation}

We also recall the definition of the first Gaussian map (cf.
\cite{wahl1}),
$$\mu_{1,L}: \Lambda^2 H^0(L) \rightarrow H^0(2L \otimes K_X),$$
as the restriction of $\tilde{\mu}_{1,L}$ to $\Lambda^2 H^0(L) \subset H^0(S,
L_S(-\Delta))$.
In local coordinates, if $s_0, s_1 \in H^0(L)$, $s_i = g_i l$, where $l$ is a
local section of $L$, we have
$$\mu_{1,L}(s_0 \wedge s_1) = (g_0 g_1' - g_1 g_0') l^2 dz.$$
 Moreover, we have
\begin{equation}
\label{ram} div(\mu_{1,L}(s_0 \wedge s_1)) = 2F + R,
\end{equation}
where $F$ is the bases locus of $|\langle s_0,s_1 \rangle| \subset
|H^0(L)|$, and $R$ is the ramification divisor of the induced
morphism. (See for example \cite{cm}, \cite{wahl2})

\begin{REM}
\label{rango4}
Recall that there is the following bijection:
$$ \{[Q] \in \proj(I_2(K_X)) \ | \ rk(Q) \leq 4\} \leftrightarrow $$
$$\{\{L, K_X-L,V, W\} \ | \ V \subset H^0(L), \ dim V = 2,
W \subset H^0(K_X - L), \ dim W = 2\}.$$ $Rank(Q) = 3$ if and only
if $2L = K_X$, and $V = W$ (see for example \cite{acgh} p.261).
\end{REM}

\begin{LEM}
\label{mu2}

If a quadric $Q$ of rank at most 4 corresponds to $\{L,K_X-L,V,
W\}$ and $V = \langle s_0, s_1 \rangle$, $W = \langle
t_0,t_1\rangle$, then
$$\mu_2(Q) = \mu_{1,L}(s_0 \wedge s_1) \mu_{1, K-L}(t_0 \wedge
t_1).$$ In particular $\mu_2(Q)\neq 0$.

\end{LEM}

\proof

By construction $Q = (s_0t_0) \otimes (s_1t_1) - (s_0t_1) \otimes
(s_1t_0) \in I_2(K_X)$. Locally  $s_i = g_i l$, where $l$ is a
local section of $L$, $t_i = h_i  l^{-1}dz$, so
$$\mu_2(Q) = - ( (g_0 h_0)' (g_1 h_1)' - (g_0h_1)' (h_0g_1)')(dz)^4=$$
$$ = (g_0 g_1'-g_1 g_0')(h_0 h_1'-h_1 h_0')(l^2 dz)((
l^{-1}dz)^2 dz) = \mu_{1,L}(s_0 \wedge s_1)\mu_{1,K-L}(t_0 \wedge
t_1).$$ \qed

\begin{REM}
Recall that by a theorem of Marc Green
(\cite{Gre}, see also \cite{acgh}, pp. 255), for a
non-hyperelliptic smooth curve of genus $g \geq 4$, $I_2$ is
generated by quadrics of rank $\leq 4$.
\end{REM}

We now make an easy linear algebra remark, which will be useful in the sequel.
\begin{REM}
\label{linalg} Let $X \subset {\proj}^n = {\proj}(V)$ be a
projective variety. Let $f: V \rightarrow W$ be a linear map,
$dim(W) = m+1$, $\bar{f}: {\proj}(V) = {\proj}^n \dashrightarrow
{\proj}(W) = {\proj}^m$ be the corresponding projection. Let $K$
be the kernel of $f$, and assume that ${\proj}(K) \cap X =
\emptyset$, i.e. ${\overline{f}}_{|X}$ is a morphism. This clearly
implies $dim(X) + dim(K) -1 \leq n-1$, or equivalently $rank(f)
\geq dim(X).$
\end{REM}

Consider the rational map
$$\overline{\mu}_2 : \proj(I_2(K_X)) \dashrightarrow
\proj(H^0(4K_X)),\
 [Q] \mapsto [\mu_2(Q)].$$

Let $\Gamma = \{[Q] \in \proj(I_2(K_X)) \ | \ rk(Q) \leq 4\}$, by
(\ref{mu2}) the restriction of the above map to $\Gamma$ is a
morphism.

\begin{PROP}
\label{rank}
For any curve of genus $g \geq 4$,
$$rk(\mu_2) \geq dim \Gamma + 1 \geq g-3.$$

\end{PROP}

\proof Since ${\overline{\mu}_2}_{|\Gamma} $ is a morphism, by
(\ref{linalg}) we know that $rk(\mu_2) \geq dim \Gamma + 1$. Let
us denote by ${\mathcal W} \subset W^1_{g-1}$ the subset of line
bundles $L \in W^1_{g-1}$ such that $h^0(L) =2$. ${\mathcal W}$ is
a non empty open subset of $W^1_{g-1}$ of dimension $\geq g-4$
(cf. e.g. \cite{acgh}). If we set ${\mathcal Y} := {\mathcal
W}/\langle \tau \rangle$, where $\tau$ is the involution mapping
$L$ to $K_X-L$, we can identify ${\mathcal Y}$ with a subset of
$\Gamma$. In fact, given a line bundle $L \in {\mathcal W}$, the
set $\{L, K-L, H^0(L), H^0(K-L)\}$ determines a quadric of rank
$\leq 4$ as we have seen in (\ref{rango4}).

Therefore, $dim(\Gamma) \geq dim({\mathcal W}) \geq g-4$, hence
$rk(\mu_2) \geq g-3.$

\qed

\section{Surjectivity}

In this section we give a class of examples of curves contained in the product of two curves for which the second Gaussian map is surjective, as Wahl does in Theorem 4.11 \cite{wahl1} for the first Gaussian map.
Other examples of curves whose second Gaussian map is surjective have been obtained by Ballico and Fontanari (\cite{bafo}) in the case of complete intersections, generalizing Wahl's result on the first Gaussian map for complete intersections (\cite{wahl3}).

Let $C_1$, $C_2$ be two smooth curves of respective genera $g_1$, $g_2$, denote by $K_i = K_{C_i}$, $i =1,2$, choose $D_i$  divisors on $C_i$ of degree $d_i$, $i=1,2$. Let $X = C_1 \times C_2$, $C \in | {p_1}^*(D_1) \otimes
{p_2}^*(D_2)|$  a smooth curve, where $p_i$ is the projection from $C_1 \times C_2$ on $C_i$,
$K_X(C) = {p_1}^*(K_1(D_1)) \otimes
{p_2}^*(K_2(D_2))$.

\begin{TEO}
 If either $g_1,g_2 \geq 2$, $d_i \geq 2g_i + 5$, $i=1,2$, or $g_1 \geq 2$, $g_2 =1$, $d_1 \geq 2g_1 + 5$, $d_2 \geq 7$, or  $g_2 =0$, $d_2 \geq 7$, $d_2(g_1 -1) > 2d_1 \geq 4g_1 + 10$, then $\mu_{2,K_C}$ is surjective for a smooth curve $C \in | {p_1}^*D_1 \otimes {p_2}^*D_2|$.

Therefore under these assumptions, for the general curve of genus $g = 1 + (g_2-1)d_1 + (g_1-1)d_2 + d_1 d_2$ the second Gaussian map is surjective.
 \end{TEO}

\proof

Denote by $I_2(K_X(C))$ the kernel of the multiplication map $$S^2H^0(K_X(C)) \rightarrow H^0(K_X^2(2C)).$$

Let $\mu^X_{2, K_X(C)}: I_2(K_X(C)) \rightarrow H^0(S^2 \Omega^1_X \otimes K_X^2(2C))$
be the second Gaussian map of the line bundle $K_X(C)$ on the surface $X$.
We have the following commutative diagram

\xymatrix{
I_2(K_X(C)) \ar[dd]^g\ar[r]^{\mu^X_{2, K_X(C)} \ \ \ \ \ \ } & H^0(S^2 \Omega^1_X \otimes K_X^2(2C)) \ar[dr]^{p_1} \\
&&H^0(S^2 \Omega^1_{X|C} \otimes K_C^2) \ar[dl]^{p_2}\\
I_2(K_C) \ar[r]_{\mu_2}& H^0(K_C^4)}

where $p_1$ is the restriction map, and $p_2$ is the map coming from the conormal extension.

We will prove that $p_1$, $p_2$ and $\mu^X_{2, K_X(C)}$ are surjective. From this we clearly obtain the surjectivity of $\mu_2$.

We want to show that $H^1(\Omega^1_{X|C} \otimes K_C^2(-C)) =0$, hence surjectivity of $p_2$ will follow.

$$H^1(\Omega^1_{X|C} \otimes K_C^2(-C)) = H^1(C, {\mathcal O}_C({p_1}^*(K_1^3(D_1)) \otimes
{p_2}^*(K_2^2(D_2))) )\oplus $$
$$ \oplus H^1(C,{\mathcal O}_C( {p_1}^*(K_1^2(D_1)) \otimes
{p_2}^*(K_2^3(D_2)))),$$ so it is sufficient to check that, under
our assumptions, ${\mathcal O}_C({p_1}^*(K_1^3(D_1)) \otimes
{p_2}^*(K_2^2(D_2)))$, and $ {\mathcal O}_C({p_1}^*(K_1^2(D_1))
\otimes {p_2}^*(K_2^3(D_2)) )$ have both degree greater than
$2g(C) -2 = d_1(2g_2 -2 + d_2) + d_2(2g_1 -2 + d_1)$.
 Let us now
consider the map $p_1$.

$$S^2 \Omega_X^1 \otimes K_X^2(C) = ({p_1}^*(K_1^4(D_1)) \otimes
{p_2}^*(K_2^2(D_2))) \oplus $$
$$ ({p_1}^*(K_1^2(D_1)) \otimes
{p_2}^*(K_2^4(D_2))) \oplus $$
$$  ({p_1}^*(K_1^3(D_1)) \otimes
{p_2}^*(K_2^3(D_2))),$$
thus by K\"unneth if either $g_i \geq 1$, $i =1,2$, or $g_2 =0$, $d_2 \geq 7$, $g_1 \geq 1$,
$ H^1(S^2 \Omega_X^1 \otimes K_X^2(C))  = 0,$ hence $p_1$ is surjective.

We want now to show that $\mu^X_{2, K_X(C)}$ is surjective.
Observe that
$$S^2H^0(K_X(C)) = (S^2H^0(K_1(D_1)) \otimes S^2H^0(K_2(D_2))) \oplus $$
$$(\Lambda^2 H^0(K_1(D_1))\otimes \Lambda^2H^0(K_2(D_2))),$$
so we have
$$I_2(K_X(C)) = \{(I_2(K_1( D_1)) \otimes S^2H^0(K_2(D_2))) + (S^2H^0(K_1(D_1)) \otimes I_2(K_2( D_2)))\} \oplus$$
$$\oplus (\Lambda^2 H^0(K_1(D_1))\otimes \Lambda^2H^0(K_2(D_2))).$$

Since
$$H^0(S^2 \Omega_X^1 \otimes K_X^2(2C)) = (H^0(C_1,K_1^4(2D_1)) \otimes
H^0(C_2, K_2^2(2D_2))  ) \oplus$$
 $$(H^0(C_1,K_1^2(2D_1)) \otimes
H^0(C_2, K_2^4(2D_2))  ) \oplus$$
$$(H^0(C_1,K_1^3(2D_1)) \otimes
H^0(C_2, K_2^3(2D_2))  ),$$
one can easily check that $\mu^X_{2, K_X(C)}: I_2(K_X(C)) \rightarrow H^0(S^2 \Omega_X^1 \otimes K_X^2(2C))$ is the sum of the three following maps:
$$
 \mu_{2,K_1(D_1)} \otimes m_2: I_2(K_1(D_1)) \otimes S^2H^0(K_2(D_2)) \rightarrow H^0(K_1^4(2D_1)) \otimes H^0(K_2^2(2D_2) ),
$$
$$
n_2 \otimes  \mu_{2,K_2(D_2)} :S^2H^0(K_1(D_1))\otimes I_2(K_2(D_2)) \rightarrow H^0(K_1^2(2D_1)) \otimes H^0(K_2^4(2D_2) ),
$$
$$
 \mu_{1, K_1(D_1)} \otimes \mu_{1,K_2(D_2)}: \Lambda^2(H^0(K_1(D_1)) \otimes \Lambda^2H^0(K_2(D_2)) \rightarrow H^0(K_1^3(2D_1)) \otimes H^0(K_2^3(2D_2) ),
$$
 where $m_2$ and $n_2$ are the multiplication maps.
 Now we apply Theorem (1.7) of \cite{bel} to the line bundles $L_i := K_i(D_i)$ on the curves $C_i$, $i =1,2$, which says that if $deg(L_i) =: l_i$ satisfies
 $2l_i \geq 3(2g_i + 2) + 2g_i -1,$
 then both $\mu_{2,L_i}$ and $\mu_{1, L_i}$ are surjective.
 So If $d_i \geq 2g_i + 5$, $\mu_{2,L_i}$ and $\mu_{1, L_i}$ are surjective, hence $\mu^X_{2, K_X(C)}$ is surjective, and this concludes the proof.
  \qed
 \begin{REM}
 The example of lowest genus of a smooth curve $C \in | {p_1}^*D_1 \otimes {p_2}^*D_2|$ with surjective second
 Gaussian map  is 71, obtained choosing $g_1 = 2$, $g_2 = 1$, $d_1 = 9$, $d_2 =7$.
\end{REM}

\section{Hyperelliptic and trigonal curves}

Assume now that $X$ is either a hyperelliptic curve of genus $\geq
3$, or a trigonal curve of genus $g \geq 4$. Let $|F|$ denote the
$g^1_2$ in the hyperelliptic case, the $g^1_3$ in the trigonal
case. Let $\phi_F: X \rightarrow {\proj}^1$ be the induced
morphism and $\nu: {\proj}^1 \hookrightarrow {\proj}^{g-1}$ be the
Veronese embedding, so that in the hyperelliptic case $\phi_K =
\nu \circ \phi_F$, where $\phi_K$ is the canonical map. Observe
that in the hyperelliptic case the hyperelliptic involution $\tau$
acts as $-Id$ on $H^0(K_X)$, so  we have an exact sequence
$$0 \rightarrow I_2(K_X) \rightarrow S^2(H^0(K_X)) \rightarrow H^0(2K_X)^+
\rightarrow 0,$$ where  $H^0(2K_X)^+$ denotes the $\tau$-invariant
part of $ H^0(2K_X)$ whose dimension is $(2g-1)$ and $I_2(K_X)$ is
the vector space of the quadrics containing the rational normal
curve.

Set $L := K_X -F$, and fix a basis $\{x,y\}$ of $H^0(F)$, and a
basis $\{t_1,...,t_r\}$ of $H^0(L)$ both in the hyperelliptic and
in the trigonal case. We have a linear map
$$\psi: \Lambda^2(H^0(L)) \rightarrow I_2, \ \ \ t_i \wedge t_j \mapsto Q_{ij}= xt_i \odot yt_j -
xt_j \odot y t_i.$$
 We recall that in both cases the linear map $\psi: \Lambda^2(H^0(L)) \rightarrow I_2$ is an isomorphism as can be easily
checked or found in \cite{am}.

\begin{LEM}
\label{hyptri2} Let $X$ be either a hyperelliptic curve of genus
$\geq 3$, or a trigonal curve of genus $g \geq 4$, let $q_1, ...,
q_l$ be the ramification points of either the $g^1_2$, or of the
$g^1_3$. Then
$$\mu_2(Q) = \mu_{1,F}(x \wedge y) \mu_{1,L}(\psi^{-1}(Q)),$$
for any quadric $Q$ of rank 4.
In particular the image of $\mu_2$ is contained in $H^0(4K_X - (q_1+...+q_l))$ and $rank(\mu_2) = rank(\mu_{1,L})$.
\end{LEM}
\proof The first statement is straightforward. So we have
$$div(\mu_2(Q)) =   div(\mu_{1,F}(x \wedge y)) +
div(\mu_{1,L}(\psi^{-1}(Q))) = $$
$$=q_1+...+q_l + div(\mu_{1,L}(\psi^{-1}(Q)).$$
Therefore $\mu_2(Q)(q_i) = 0 \ \forall i =1,...,l.$\\
\qed

\begin{PROP}
\label{hyp} Let $X$ be a hyperelliptic curve of genus $g \geq 3$,
then the rank of  $\mu_2$ is $2g -5$.
\end{PROP}

\proof Given a hyperelliptic curve of genus $g$ with equation $y^2
= f(x)$, where $f$ has degree $2g+2$ and only simple roots, a
basis of $H^0(K_X)$ is  $\{\omega_i =x^i \frac{dx}{y} \ |
\ 0\leq i \leq g-1 \}$. Let $|F|$ be the $g^1_2$ on $X$ and we
assume $F = \phi_F^{-1}(0) =: p_1 + p_2$.

Set $L = K_X - F = K_X - p_1 -p_2$, $H^0(L) \subset H^0(K_X)$ and
let $\mu_{1,L}: \Lambda^2 H^0(L)  = \Lambda^2 H^0(K_X -p_1-p_2)
\rightarrow H^0(2L + K_X)= H^0(3K_X - 2p_1 -2p_2)$ be the first
Gaussian map of $L$, then $$\mu_{1,L} =
{\mu_{1,K}}_{|{\Lambda}^2H^0(K_X-p_1-p_2)}.$$ By (\ref{hyptri2}),
the rank of $\mu_2$ is equal to the rank of $\mu_{1,L}$. As it is
shown in \cite{cm},
$$\mu_{1,K}(\omega_i \wedge \omega_j) =(i-j)\frac{x^{i+j-1}}{y^2} (dx)^3, \ 0 \leq i<j \leq
g-1,$$ so there are exactly $2g -3$ distinct powers of $x$.

Then clearly a basis of $H^0(K_X - p_1 - p_2)$ is given by $\{x^i \frac{dx}{y},
\ i>0\}$. We want to compute the dimension of the span of
$\{\mu_{1,K} (\omega_i \wedge \omega_j), \ 0<i<j \leq g-1\}$.
Observe that  $l: = i+j -1 =0$ if and only if $i = 0$, $j = 1$; $l
=1$ if and only if $i = 0$, $j =2$. But if $l \geq 3$, $l=i + j
-1$ also for some $i,j>0$.

Hence  $rank(\mu_{1,L})=rank(\mu_{1,K})-2=2g-5$.  \qed
\\

Assume now that $X$ is non-hyperelliptic trigonal curve of genus
$g \geq 4$. Let $|F|$ be the $g^1_3$ on $X$, assume $F = p_1 + p_2
+ p_3$, $ p_i \in X$. Let us denote by $L = K_X -F = K_X
-p_1-p_2-p_3$, $deg(L) = 2g-5$, $h^0(L) = g-2$. So $H^0(L) \subset
H^0(K_X)$ and $\mu_{1,L} =
{\mu_{1,K}}_{|{\Lambda}^2H^0(K_X-p_1-p_2-p_3)}.$
 In \cite{cm} it is proven
that for the general trigonal curve of genus $g \geq 4$,
$dim(coker(\mu_{1,K})) = g +5$, moreover specific examples of
trigonal curves (whose genera are all equal to 1 modulo 3) such
that the corank of $\mu_{1,K}$ is $g+5$ are exhibited. Using
results of \cite{md}, in \cite{br} Brawner  proves that
$dim(coker(\mu_{1,K})) = g +5$ for any trigonal curve of genus $g
\geq 4$.

We shall now compute the rank of $\mu_2$ for trigonal curves. By
(\ref{hyptri2}) it suffices to compute $rank(\mu_{1,L})$ and we
will do it following the computation done in \cite{md} and
\cite{br} for $\mu_{1,K}$.

Recall that a canonically embedded trigonal curve of genus $g$ lies on a rational normal scroll $S_{k,l}$, where $k \leq l$, $l +k = g-2$ and $k$ is the Maroni invariant, which is bounded by
\begin{equation}
\label{maroni}
\frac{g-4}{3} \leq k \leq \frac{g-2}{2},
\end{equation}
(cf. \cite{ma}).

The surface $S_{k,l}$ is isomorphic to ${\bf F}_n$, with $n =
l-k$. Let us denote by $H$ the hyperplane section, by $R$ the
fibre of the ruling, and set $B \equiv H - lR$.

We have
$$H^2 = g-2, \ B^2 = -n,$$
$$C \equiv 3H -(g-4)R,$$
$$K_S \equiv -2H + (g-4)R,$$
hence
 $$K_S + C -R \equiv H-R \equiv B + (l-1)R,$$
and
$$(K_S + C -R)_{|C} \equiv L. $$

\begin{TEO}
\label{trigonali} For any trigonal curve $C$ of genus $g \geq 8$,
the rank of $\mu_2$ is $4g -18$. Hence the general curve of genus
$g \geq 8$, $\mu_2$ has rank greater or equal to $4g-18$.
\end{TEO}

\proof

As in \cite{md}, (2.1)  we have the following commutative diagram
involving the first Gaussian map $\mu_{1, H-R}^{S}$ for the scroll
$S := S_{k,l}$.
$$
\begin{array}{ccc}
  \Lambda^2H^0(S, {\mathcal O}_S(K_S + C -R))& \stackrel{\mu_{1, H-R}^{S}}\longrightarrow& H^0(S, \Omega^1_S(2K_S + 2C -2R)) \\
  Res\downarrow&  & \gamma'\downarrow \\
  \Lambda^2H^0(C, K_C -F)) & \stackrel{\mu_{1, L}}\longrightarrow& H^0(C, 3K_C -2F) \\
\end{array}
$$

We will prove that the map $\mu_{1, H-R}^{S}$ is surjective,
$\gamma'$ is injective and $Res$ is surjective. This implies that
$rank(\mu_{1,L}) = h^0(S, \Omega^1_S(2K_S + 2C -2R)) = h^0(S,
\Omega^1_S(2B + 2(l-1)R))$.

Observe that by the bound (\ref{maroni}) of the Maroni invariant,
$k \geq 2$ for $g \geq 8$,  hence  the hypothesis of corollary
(3.3.2) of \cite{md} are satisfied and thus  $h^0(S, \Omega^1_S(2B
+ 2(l-1)R)) = 4g-18$. In fact corollary (3.3.2) of \cite{md}
asserts that  $h^0(S, \Omega^1_S(rB + sR)) = 2rs -nr^2 -2$, if $r
\geq 1$, and $s \geq nr+2$.

The surjectivity of  $\mu_{1, H-R}^{S}$ follows by theorem (4.5)
of \cite{md} that says that $\mu_{1, rB +sR}^{S}$ is surjective if
$r \geq 0$ and $s \geq nr +1$.

In \cite{br}, (3.4) it is proven that the map
$$\gamma: H^0(S, \Omega^1_S(2H)) \rightarrow H^0(C, 3K_C)$$
is injective. Since $\gamma'$ is the restriction  of $\gamma$ to
$H^0(S, \Omega^1_S(2H-2R))$, also $\gamma'$ is injective.

We finally show that the restriction map
$$H^0(S, {\mathcal O}_S(H-R)) \rightarrow  H^0(C, L)$$
is surjective.
Consider the exact sequence
$$0 \rightarrow {\mathcal O}_S(H-R-C) \rightarrow {\mathcal O}_S(H-R) \rightarrow {\mathcal O}_C(H-R) \rightarrow 0.$$
An easy computation on the scroll shows that $H^1(S, {\mathcal
O}_S(H-R-C)) = H^1(S, {\mathcal O}_S(-2H + (g-5)R)) =0$, proving
our assertion.
\qed\\

\section{Injectivity for low genus}
We will now exhibit some examples of computations of the rank of
$\mu_2$ for genus $\leq 7$ from which it will follow that $\mu_2$
is injective for the general curve of genus $\leq 6$. Note that if
$g(X) =4$, $I_2(K_X)$ has dimension 1, so $\mu_2$ is injective.

\begin{PROP}
\label{tri5} For any trigonal curve $X$ of genus 5, $\mu_2$ is
injective. Hence for the general curve of genus 5 $\mu_2$ is
injective.
\end{PROP}
\proof For a curve of genus 5 the dimension of $I_2(K_X)$ is 3.
Let us assume that $X$ is trigonal. Then there exists a line
bundle $L$ on $X$ such that $h^0(L) = 2$, $deg(L)=3$, so that
$h^0(K-L) = 3$. Let $Gr(2,H^0(K-L))$ be the Grassmannian of the 2
dimensional subspaces in $H^0(K-L)$. To any $W \in Gr(2,H^0(K-L))$
we associate the quadric $Q_W $ of rank $4$ corresponding to the
set $\{L, K-L, H^0(L), W\}$ as in (\ref{mu2}). Thus we have a
morphism
$$Gr(2, H^0(K-L)) \rightarrow \proj(H^0(4K)), \ \ \ W \mapsto \overline{\mu}_2(Q_W). $$ So by (\ref{linalg}) we have

$rank(\mu_2) \geq dim( Gr(2,H^0(K-L)) + 1 = 3$
\qed\\

\begin{TEO}
\label{quintic} Let $X$ be a smooth plane quintic, then the map
$\mu_2$ is injective and its image has no base points. Hence for
the general curve of genus 6 $\mu_2$ is injective.
\end{TEO}
\proof $X$ has genus 6, so the dimension of $I_2(K_X)$ is 6. We
will find 6 quadrics of rank at most 4 such that their images
under $\mu_2$ are linearly independent. Observe that $K_X \equiv
{\mathcal O}_X(2)$, hence $L := {\mathcal O}_X(1)$ is such that
$2L \equiv K_X$. Let $q_1, q_2, q_3$ be distinct points of $X$ in
general position such that the tangent line $r_i$ of $X$ at $q_i$
is a simple tangent, $i=1,2,3$. Assume furthermore that $P_{12} :=
r_1 \cap r_2$, $P_{13} := r_1 \cap r_3$, $P_{23} := r_2 \cap r_3$
are in general position and do not lie on $X$. Denote by
$\pi_{12}$, $\pi_{13}$, $\pi_{23}$ the respective projections $X
\rightarrow {\proj}^1$ from the points $P_{12}$, $P_{13}$,
$P_{23}$. These 3 projections $\pi_{ij}$ correspond to 3 pencils
$V_{ij} \subset H^0(L)$ with $2L \equiv K$. Let $R_{ij}$ (for $1
\leq i<j \leq 3$) be the ramification divisor of $\pi_{ij}$, then
we have
$$R_{ij} = q_i + q_j + A_{ij}.$$
Observe that by construction $q_3 \not\in A_{12}$, otherwise
$P_{12}$ would lie on  $r_3$ and we would have $P_{12} = P_{13} =
P_{23}$. Analogously $q_1 \not \in A_{23}$, $q_2 \not \in A_{13}$.
Observe that, since $r_2$ and $r_3$ are simple tangents, $q_2,q_3
\not\in A_{23}$, hence there must exist a point $q_4 \in A_{23}$,
which is different from $q_2$ and $q_3$. Notice that $q_4 \not\in
R_{12} \cup R_{13}$. In fact, by construction $P_{23} \in r_4$
($r_4$ is the tangent line of $X$ at $q_4$), hence $P_{23} = r_2
\cap r_3 = r_4 \cap r_2 = r_4 \cap r_3$. So if $q_4 \in R_{12}$,
then $P_{12} \in r_4$, thus $P_{12} = r_4 \cap r_2 = P_{23}$, a
contradiction. Analogously, if $q_4 \in R_{13}$ we get $P_{12} =
P_{23} = P_{13}$, which is impossible.

Define now the 6 quadrics in $\Gamma$ by the following 6 sets as
in (\ref{mu2}):
$$Q_{ij,kl} \longleftrightarrow \{L,L=K-L,V_{ij}, V_{kl} \}$$
$1 \leq i<j \leq 3, 1 \leq k<l \leq 3.$

We have by (\ref{mu2})
$$div(\mu_2(Q_{12,12})) = 2R_{12}= 2q_1 + 2q_2 + 2A_{12},$$
$$div(\mu_2(Q_{12,13})) = R_{12} + R_{13} = 2q_1 + q_2 + q_3 +
A_{12} + A_{13},$$
$$div(\mu_2(Q_{12,23})) = R_{12} + R_{23} = q_1 + 2q_2 + q_3 +
A_{12} + A_{23},$$
$$div(\mu_2(Q_{13,13})) = 2R_{13} = 2q_1 + 2q_3 + 2A_{13},$$
$$div(\mu_2(Q_{13,23})) = R_{13} + R_{23} = q_1 + q_2 + 2q_3 +
A_{13} + A_{23},$$
$$div(\mu_2(Q_{23,23})) = 2R_{23} = 2q_2 + 2q_3 + 2A_{23}.$$

Assume now that there exists a linear combination $\sum
\lambda_{ij,kl}\mu_2(Q_{ij,kl}) =0,$ then by evaluating in $q_1$
we get $\lambda_{23,23} \mu_2(Q_{23,23})(q_1) =0$, thus
$\lambda_{23,23}=0$, since $q_1 \not \in R_{23}$. Evaluation in
$q_2$ yields $\lambda_{13,13} \mu_2(Q_{13,13})(q_2) =0$, thus
$\lambda_{13,13}=0$, since $q_2 \not \in R_{13}$. By evaluating in
$q_3$ we obtain $\lambda_{12,12} \mu_2(Q_{12,12})(q_3) =0$, thus
$\lambda_{12,12}=0$, since $q_3 \not \in R_{12}$. We now evaluate
in $q_4 \in A_{23}$ and we have $\lambda_{12,13}
\mu_2(Q_{12,13})(q_4) =0$. This implies $\lambda_{12,13} =0$,
since otherwise we would have $q_4 \in A_{12} + A_{13}$, which is
impossible. Finally we must have $\lambda_{12,23} =
\lambda_{13,23} =0$, otherwise we would have $R_{12} = R_{13}$, a
contradiction. This proves injectivity of $\mu_2$.

Notice that if $P$ is a base point of the image of $\mu_2$, the
three projections $\pi_{ij}$ must have a common ramification
point, which is impossible by construction.

 \qed

If $X$ is a trigonal curve of genus 7, the argument of
(\ref{tri5}) yields $rank(\mu_2)\geq dimGr(2,H^0(K-L))+1=7$. We
will now exhibit an example of a trigonal curve of genus 7 such
that $rank(\mu_2) = rank(\mu_{1,L}) = 9.$ Our example is the
cyclic covering of ${\proj}^1$ constructed in \cite{cm}, whose
affine equation is
$$y^3= x^{9} -1.$$
In \cite{cm} the map $\mu_{1,K}$ is explicitly computed on the
elements ${\sigma}_{ij} = x^i y^j \frac{dx}{y^2}$, $0\leq j \leq
1$, $0 \leq i \leq 3(2-j) -2$, which form a basis of $H^0(K_X)$.
Namely it is shown that the image of $\mu_{1,K}$ (of dimension 18)
is spanned by the following 3 types of elements:
\begin{equation}
\label{type1} \mu_{1,K}(\sigma_{i0} \wedge \sigma_{k0}) =
[(k-i)x^{i+k-1} y^{-4}] \cdot (dx)^3, \ 0 \leq i <k \leq 4,
\end{equation}
\begin{multline}
\label{type2} \mu_{1,K}(\sigma_{i0} \wedge \sigma_{k1}) =
[(k-i)x^{i+k-1} y^{-3} + 3x^{i+k+8} y^{-6}] \cdot
(dx)^3, \\
0 \leq i \leq 4, \ 0 \leq k \leq 1,
\end{multline}
\begin{equation}
\label{type3} \mu_{1,K}(\sigma_{01} \wedge \sigma_{11}) = y^{-2}
\cdot (dx)^3.
\end{equation}
Now the $g^1_3$ on our curve is the linear system $|F| = |p_1 +
p_2 + p_3|$ where $p_i$ can be chosen to be the points $(0, y_i)$,
with $y_i^3 = -1$. So, if $L = K_X - p_1-p_2 -p_3$  we can
identify $H^0(L)$ with the subspace of $\langle \sigma_{ij}
\rangle$ generated by the elements $\sigma_{ij}$, where $i >0$.
Since $\mu_{1,L} = {\mu_{1,K}}_{|\Lambda^2 H^0(L)}$, one has to
compute the dimension of $ \langle \mu_{1,K}(\sigma_{ij} \wedge
\sigma_{kl} ) \ | \ i, k>0 \rangle $ which turns out to be 9.

\section{Base points} We will now show global generation of the
image of $\mu_2$ for curves which are not hyperelliptic nor
trigonal.

\begin{TEO}
\label{puntobase} Assume that $X$ is smooth curve of genus $g \geq
5$, which is non-hyperelliptic and non-trigonal. Then for any $P
\in X$ there exists a quadric $Q \in I_2$ such that $\mu_2(Q)(P)
\neq 0$. Equivalently $Im(\mu_2)\cap H^0(4K_X-P) \neq Im(\mu_2)$,
$\forall P \in X.$
\end{TEO}
\proof We will show that for any $P \in X$ there exists a quadric
$Q$ of rank 4 such that $\mu_2(Q)(P) \neq 0$. We recall that any
component of the space $ W^1_{g-1}(X)$ has dimension greater or
equal to $g-4$, and in \cite{ckm} lemma (2.1.1) it is proven that
if $X$ is non-hyperelliptic, non-trigonal and not isomorphic to a
smooth plane quintic, there exists a line bundle $L \in W^1_{g-1}$
such that both $|L|$ and $|K_X -L|$ are base point free. If $X$ is
a plane quintic, by (\ref{quintic}) we know that $\mu_2$ has no
base points. So we can assume that there exists a non empty
irreducible open subset ${\mathcal V}$ in $ W^1_{g-1}(X)$ which
consists of line bundles $L$ such that $h^0(L) = 2$, $L \not\equiv
K_X-L$, both $|L|$ and $|K_X -L|$ are base point free. So the
condition $\mu_2(Q)(P) = 0$ for the quadric associated to $|L|$
and $|K_X -L|$ says that either $P$ is a ramification point for
the morphism $\phi_{|L|}:X \rightarrow \proj^1$ or for the
morphism $\phi_{|K_X-L|}:X \rightarrow \proj^1$.

We claim that there exists $L \in  {\mathcal V}$ such that $P$ is
at most a simple ramification point for both $\phi_{|L|}$ and
$\phi_{|K_X-L|}$, that is $h^0(L-3P) =0$ and $h^0(K-L-3P)=0$. In
fact, assume that for all $L \in  {\mathcal V}$ either $h^0(L-3P)
\geq1$, or $h^0(K-L-3P)\geq 1$. Consider the maps $F_1:
Sym^{g-4}(X) \rightarrow Pic^{g-1}(X),$ $F_1(D) = D + 3P$, $F_2:
Sym^{g-4}(X) \rightarrow Pic^{g-1}(X),$ $F_2(D) = K_X - D - 3P$.
Then  ${\mathcal V}$ is contained in $Im(F_1) \cup Im(F_2)$, and
since they have the same dimension and  ${\mathcal V}$ is
irreducible, $\overline{\mathcal V} = Im(F_1)$ or
$\overline{\mathcal V} = Im(F_2)$. This means that $\forall
x_1,...,x_{g-4} \in X$, $h^0(x_1+...+x_{g-4} + 3P) = h^0(K_X
-x_1-...-x_{g-4} - 3P) \geq 2$, which is absurd, since
$h^0(K_X-3P) =g-3$, because $X$ is non hyperelliptic and non
trigonal.

So assume that $P$ is a simple ramification point for $\phi :=
\phi_{|L|}: X \rightarrow \proj^1$, let $R$ be its ramification
divisor and consider the exact sequence
$$0 \rightarrow T_X \stackrel{\phi_*} \rightarrow \phi^* T_{\proj^1}
\rightarrow {\mathcal N}_{\phi} \rightarrow 0.$$
Notice that $ \phi^* T_{\proj^1} = T_X(R)$, $ {\mathcal N}_{\phi} = T_X(R)_{|R}$ and $\phi_*$ is the inclusion map of $T_X$ in $T_X(R)$, so that our exact sequence is
 $$0 \rightarrow T_X \rightarrow T_X(R)
\rightarrow  T_X(R)_{|R} \rightarrow 0.$$

We have a cohomology exact sequence
$$0 \rightarrow H^0(T_X(R))
\stackrel{i} \rightarrow H^0( T_X(R)_{|R})
\stackrel{\beta}\rightarrow H^1(T_X) \rightarrow 0.$$ We recall
(see e.g.  \cite{ac}) that $H^0(T_X(R)) =  H^0(\phi^*
T_{\proj^1})$ parametrises the infinitesimal deformations of the
morphism $\phi$ which do not move $X$ and it contains the three
dimensional subspace $\phi^*(H^0(T_{\proj^1}))$. Note that since
$H^1(T_X(R))=0$ any first order deformation extends.

The strategy of the proof is to exhibit an infinitesimal
deformation $\rho \in H^0(T_X(R))$ of $\phi$ such that $i(\rho)
\not \in H^0(T_X(R-P)_{|R-P}) \subset H^0( T_X(R)_{|R})$ and
$\rho$ is not contained in $\phi^*(H^0(T_{\proj^1}))$. The
assumption that the ramification in $P$ is simple implies that
$\rho$ extends to a deformation of $\phi$ such that the new map is
not ramified in $P$ anymore.

Denote by $\xi_{P} \in H^1(T_X)$ a Schiffer variation in $P$, that
is, by definition a generator of the subspace
$Im(H^0(T_X(P))_{|P}) \subset H^1( T_X)$. Since $P$ is a
ramification point, $\xi_P\in \beta(H^0(T_X(R)_{|R}))$.  Notice
that if we prove that $\beta(H^0(T_X(R-P)_{|R-P}))\subset
\beta(H^0(T_X(R)_{|R}))$ generates $H^1(T_X)$, then there exists
an element $x \in H^0(T_X(R-P)_{|R-P})$ such that $\xi_P =
\beta(x)$. So if $c_P \in H^0(T_X(P)_{|P}) \subset
H^0(T_X(R)_{|R})$ is such that $\beta(c_P) = \xi_P$, the element
$\eta = c_{P} - x \in H^0(T_X(R)_{|R})$ maps to zero in
$H^1(T_X)$, thus there exists an element $\theta \in H^0(T_X(R))$
such that $\eta = i(\theta)$ and $\eta \not \in
H^0(T_X(R-P)_{|R-P})$ since the coefficient of $c_{P}$ in $\eta$
is non zero.

Therefore we want to prove that there exists $L \in {\mathcal V}$
with simple ramification in $P$ such that
$\beta(H^0(T_X(R-P)_{|R-P}))$ generates $H^1(T_X)$. Assume to the
contrary that for any $L \in {\mathcal V}$,
$\beta(H^0(T_X(R-P)_{|R-P}))$ lies on a hyperplane in $H^1(T_X)$,
i.e. there exists an element $\omega \in H^0(2K_X)$ such that
$\omega(P_i) = 0$ for all $i \geq 2$, and $ord_{P_i} \omega =
n_i$, where $\sum_{i\geq 2} n_i P_i=R-P$. Then we have
$$2K_X \equiv div(\omega)= P + \sum_{i \geq 2} n_iP_i -P + q \equiv K_X + 2L
-P +q,$$ for some $q \in X$. So for any $L$ there exists a point
$q \in X$ such that
$$2L \equiv K_X + P -q.$$
Since $L$ varies in an open subset of $W^1_{g-1}$, which has
dimension at least $g-4$, while $q$ varies in $X$, if $g \geq 6$
this cannot hold for all $L \in {\mathcal V}$.

If $g =5$ we still get a contradiction noting that the
multiplication by 2 in the jacobian restricts to a connected
topological covering $\tilde{X}$ of the curve $X$ of degree
$2^{10}$ corresponding to the surjective homomorphism
$\pi_1(X)\rightarrow H_1(X,\Z/2\Z)$. Hence $\tilde{X}$ can't
coincide with a component of $W^1_4$ which is a 2 to 1 covering of
a quintic plane curve ramified along at most 10 points (cf.
\cite{acgh} p.270).

Finally we show that we can choose the deformation outside
$\phi^*(H^0(T_{\proj^1}))$. Set $W = i(H^0(T_X(R)) \cap
H^0(T_X(R-P)_{|R-P})$. We have just proven that, for $L$ general,
$i(H^0(T_X(R)) \not \subset H^0(T_X(R-P)_{|R-P})$, so $dim(W) =
h^0(T_X(R))  -1 =
  g-2 \geq 3$, since $g \geq 5$. Thus if $\{e_1, e_2, e_3\}$ are three linearly
  independent elements in $W$, the four elements $\{i(\theta) =\eta, \eta + e_1,
  \eta+e_2, \eta +e_3\}$ are linearly independent, since $\eta$ is not
  contained in $H^0(T_X(R-P)_{|R-P})$ therefore there exists a deformation $\rho$ such that $i(\rho)  \in \{i(\theta)
=\eta, \eta +
  e_1, \eta+e_2, \eta +e_3\ \} \subset i(H^0(T_X(R)))$ which does
  not belong to the 3 dimensional subspace
  $i(\phi^*(H^0(T_{\proj^1})))$. Hence $\rho$ is the deformation
  we are looking for.

We have proven that if $L$ does not belong to the curve $\gamma$
given by the equation $2L \equiv K_X + P -q,$ with $q \in X$, we
can deform $L$ in such a way that $P$ is not a ramification point
anymore.  Analogously if $L_1 := K_X -L$ does not belong to the
curve $\gamma$ we can deform $L_1 = K_X - L$ in such a way that
$P$ is not a ramification point of the corresponding morphism
anymore. So if we take $L \in {\mathcal V} - {\gamma} -
\iota^{-1}(\gamma)$, with $\iota: W^1_{g-1} \rightarrow W^1_{g-1}$
the involution sending $L$ to $K_X -L$, we find deformations of
$L$ (then also of $K_X-L$) such that if $L'$ is the deformed line
bundle, $P$ is
neither a ramification point of $\phi_{|L'|}$ nor of $\phi_{|K_X-L'|}$.  \\

\qed\\


\begin{thebibliography} {99}

\bibitem{am} Andreotti,~A., Mayer,~A., On period relations for abelian integrals on algebraic curves, {\em
    Ann. Scuola Norm. Sup. Pisa} {\bf 21}
(1967), no. 2, 189--238.

\bibitem{ac} Arbarello,~E., Cornalba,~M., Su una congettura di Petri, {\em
    Comment. Math. Helvetici} {\bf 56}
(1981), no. 1, 1--38.

\bibitem{acgh} Arbarello,~E., Cornalba,~M., Griffiths,~P., Harris, ~J. {\em
Geometry of algebraic curves, Vol. I}, Grundlehren der Mathematischen Wissenschaften, 267. Springer-Verlag, New York, 1985.

\bibitem{bafo} Ballico, ~E., Fontanari, ~C.,  {\em On the surjectivity of higher Gaussian maps for complete intersection curves}.  Ricerche Mat.  53  (2004),  no. 1, 79--85 (2005).

\bibitem{bel} Bertram,~A., Ein,~L., Lazarsfeld,~R. Surjectivity of Gaussian Maps for Line Bundles of Large Degree on Curves in {\em Algebraic Geometry}, Chicago 1989,
Lecture Notes in Mathematics, 1479. Springer, Berlin, (1991),
15--25.

\bibitem{br} Brawner,~J., The Gaussian-Wahl map for trigonal curves,  {\em Proc. Amer. Math. Soc.} {\bf 123}, no.\ 5,
(1995), 1357--1361.

\bibitem{chm} Ciliberto,~C., Harris,~J., Miranda, ~R. On the surjectivity of the Wahl map, {\em Duke Math. Jour.}, {\bf 57} , (1988), 829--858.

 \bibitem{cm} Ciliberto,~C., Miranda, ~R., Gaussian maps for certain families
  of canonical curves, {\em Complex projective geometry} (Trieste,
  1989/Bergen, 1989), London Math. Soc. Lecture Note Ser. 179, Cambridge
  Univ. Press, Cambridge, (1992) 106--127.

 \bibitem{cm1} Ciliberto,~C., Miranda, ~R., Gaussian map for canonical curves of low genus, {\em Duke Math. Jour.}, {\bf 61} no. 2, (1990), 417--443.

\bibitem{cfcl} Colombo,~E., Frediani,~P., Siegel metric and curvature of the moduli space of curves, preprint.

\bibitem{cpt} Colombo,~E., Pirola,~G.P., Tortora,~A., Hodge-Gaussian
maps, {\em Ann. Scuola Normale Sup. Pisa Cl. Sci. (4)} {\bf 30}
(2001), no. 1, 125--146.


\bibitem{ckm} Coppens,~M., Keem,~C., Martens,~G., Primitive linear series on
  curves, {\em Manuscripta Math.} {\bf 77}
(1992), 237--264.

\bibitem{md} Duflot,~J., Miranda,~R.,  The Gaussian map for rational ruled surfaces, {\em Trans. Amer. Math. Soc.} {\bf 330}
(1992), 447--459.


\bibitem{Gre} Green,~M.~L., Quadrics of rank four in the ideal of a canonical
  curve, {\em  Invent. Math.} {\bf 75}
(1984), n.1, 85--104.

\bibitem{green} Green,~M.~L., Infinitesimal methods in Hodge theory, in {\em Algebraic Cycles and Hodge Theory}, Torino 1993,
Lecture Notes in Mathematics, 1594. Springer, Berlin, (1994),
1--92.

\bibitem{Gr} Griffiths,~P.~A., Infinitesimal variations of Hodge structures
  (III): determinantal varieties and the infinitesimal invariant of normal
  functions, {\em Comp. Math.} {\bf 50}
(1983), 267--324.


\bibitem{ma} Maroni,~A., Le serie lineari sulle curve trigonali, {\em Ann. Mat. Pura Appl.} {\bf 25} (1946), 341--354.

\bibitem{mm} Mori,~S., Mukai,~S., The uniruledness of the moduli space of curves of genus 11,
in {\em Algebraic Geometry Proceedings Tokyo, Kyoto}, 1982,
Springer LNM {\bf 1016} (1983), 334--353.

\bibitem{ser} Sernesi,~E., Moduli of rational fibrations, arXiv:math/0702865v2.


\bibitem{voi} Voisin,~C., Sur l'application de Wahl des courbes satisfaisant la condition de Brill-Noether-Petri, {\em Acta Math.} {\bf 168} (1992), 249--272.


\bibitem{wahl1} Wahl,~J., Gaussian maps on algebraic curves, {\em
J. Diff. Geom.} {\bf 32} (1990), no. 1, 77--98.

\bibitem{wahl2} Wahl,~J., Introduction to Gaussian maps on an algebraic curve,
  {\em Complex projective geometry} (Trieste, 1989/Bergen, 1989),
London Math. Soc. Lecture Note Ser., 179, Cambridge Univ. Press, Cambridge, (1992), 304--323.

\bibitem{wahl3} Wahl,~J., The Jacobian algebra of a graded Gorenstein singularity, {\em
Duke Math. Jour.} {\bf 55} (1987), 843--871.



\end{thebibliography}
\end{document}